%% file: complexSCM.tex
\documentclass[journal]{IEEEtran}
\usepackage{amssymb,amsmath,amsthm}

\newcommand{\bo}[1]{\mathbf{#1}} 
\newcommand{\bom}[1]{\boldsymbol{#1}}   

\newcommand{\fC}{\mathbb{C}}            
\newcommand{\fR}{\mathbb{R}}            
\newcommand{\im}{\jmath}

\newcommand{\Hm}[1]{\mathcal H^{#1}}
\newcommand{\PSDHm}[1]{\mathcal H_{+}^{#1}}
\newcommand{\PDHm}[1]{\mathcal H_{++}^{#1}}

\newcommand{\MSE}{\mathrm{MSE}}
\newcommand{\NMSE}{\mathrm{NMSE}}

\newcommand{\ndim}{n}                   
\newcommand{\pdim}{p}                   

\renewcommand\a{\mathbf{a}}

\renewcommand\u{\mathbf{u}}
\newcommand\z{\mathbf{z}}
\newcommand\e{\mathbf{e}}
\newcommand\q{\mathbf{q}}

\newcommand\x{\mathbf{x}}

\newcommand\M{{\boldsymbol \Sigma}}     
\newcommand\bmu{{\boldsymbol \mu}}      

\renewcommand\S{{\mathbf S}}            
\newcommand\I{{\mathbf I}}              

\newcommand{\ka}{\kappa}                
\newcommand{\kurt}{\mathrm{kurt}}       
\newcommand{\commat}[1]{\bo K_{#1,#1}}  
\newcommand{\ve}{\mathrm{vec}}          
\newcommand\veco[1]{\mathrm{vec}\left({#1}\right)} 

\newcommand\Fro{{\mathrm{F}}}           
\DeclareMathOperator{\E}{\mathbb{E}}
\DeclareMathOperator{\var}{\mathrm{var}}
\DeclareMathOperator{\pvar}{\mathrm{pvar}}
\DeclareMathOperator{\cov}{\mathrm{cov}}
\DeclareMathOperator{\pcov}{\mathrm{pcov}}
\DeclareMathOperator{\tr}{tr}

\newcommand{\hop}{\mathsf{H}}           

\usepackage{tikz}
\usepackage{pgfplots}
\pgfplotsset{compat=1.13}

\usepackage{color}

\theoremstyle{remark}
\newcounter{ctheorem}
\newtheorem{theorem}[ctheorem]{Theorem}

\newcommand{{\papertitle}}{On the variability of the sample covariance matrix
under complex elliptical distributions}
\newcommand{{\keywords}}{Affine equivariant statistics, sample covariance
matrix, complex-valued elliptical distribution, complex-valued multivariate
normal distribution}

\usepackage[pdftex,
pdfauthor={Elias Raninen, Esa Ollila, and David E. Tyler},
pdftitle={\papertitle},
pdfkeywords={\keywords},
pdfcreator={pdflatex}]{hyperref}

\title{\papertitle}
\author{Elias~Raninen,~\IEEEmembership{Student~Member,~IEEE},~Esa~Ollila,~\IEEEmembership{Senior~Member,~IEEE}
and~David~E.~Tyler \thanks{Elias Raninen and Esa Ollila are with the Department
of Signal Processing and Acoustics, Aalto University, FI-00076 Aalto, Finland.
David E. Tyler is with the Department of Statistics, Rutgers - The State
University of New Jersey, Piscataway, NJ 08854, USA. The work of David E. Tyler
was supported in part by the National Science Foundation under Grant
DMS-1812198.}
}

\begin{document}
\maketitle

\begin{abstract}
    We derive the form of the variance-covariance matrix for any affine
    equivariant matrix-valued statistics when sampling from complex elliptical
    distributions. We then use this result to derive the variance-covariance
    matrix of the sample covariance matrix (SCM) as well as its theoretical mean
    squared error (MSE) when finite fourth-order moments exist. Finally,
    illustrative examples of the formulas are presented.
\end{abstract}

\begin{IEEEkeywords}
sample covariance matrix, sample variation, mean squared error, complex Gaussian
distribution, complex elliptically symmetric distribution
\end{IEEEkeywords}

\section{Introduction}
Suppose we observe independent and identically distributed (i.i.d.)
complex-valued $p$-variate random vectors $\x_1,\ldots,\x_n \subset \fC^p$ with
mean $\bmu = \E[\x_i]$ and positive definite covariance matrix $\M = \E[(\x_i -
\bmu)(\x_i - \bmu)^\hop]$. The (unbiased) estimators of $\M$ and $\bmu$ are the
sample covariance matrix (SCM) and the sample mean defined by
\begin{equation}
    \S = \frac{1}{n-1}\sum_{i=1}^n (\x_i - \bar \x)(\x_i - \bar \x)^\hop
    ~\text{and}~
    \bar \x = \frac{1}{n}\sum_{i=1}^n \x_i.
\end{equation}
The SCM is an integral part of many statistical signal processing methods such
as adaptive filtering (Wiener and Kalman filters), spectral estimation and
array processing (MUSIC algorithm, Capon
beamformer)~\cite{caponHighresolutionFrequencywavenumberSpectrum1969,
manolakisStatisticalAdaptiveSignal2005a}, and adaptive radar detectors
\cite{kellyAdaptiveDetectionAlgorithm1986, robeyCFARAdaptiveMatched1992,
krautAdaptiveSubspaceDetectors2001}.

In signal processing applications, a typical assumption would be to assume that
the data $\x_1,\ldots,\x_n$ follow a (circular) complex multivariate normal
(MVN) distribution \cite{goodmanStatisticalAnalysisBased1963}, denoted by $\fC
\mathcal N(\bom \mu, \M)$. However, a more general assumption would be to assume
a Complex Elliptically Symmetric (CES)
\cite{krishnaiahComplexEllipticallySymmetric1986,
ollilaComplexEllipticallySymmetric2012a} distribution, which is a family of
distributions including the MVN distribution as well as
heavier-tailed distributions such as the $t$-, $K$-, and the inverse Gaussian
distribution that are commonly used in radar and array signal
processing applications as special
cases~\cite{conteModellingSimulationNonRayleigh1991,
billingsleyStatisticalAnalysesMeasured1999,
ollilaComplexEllipticallySymmetric2012a,
ollilaCompoundGaussianClutterModelling2012}.

In the paper, we study the complex-valued (unbiased) SCM for which we derive the
variance-covariance matrix as well as the theoretical mean squared error (MSE)
when sampling from CES distributions. We also provide a general expression for
the variance-covariance matrix of any affine equivariant matrix-valued statistic
(of which the SCM is a particular case). The results regarding the SCM extend
the results in~\cite{ollilaOptimalShrinkageCovariance2019a} to the
complex-valued case, where the variance-covariance matrix and MSE of the SCM
were derived for real-valued elliptical distributions.

The structure of the paper is as follows. Section~\ref{sec:CES} introduces CES
distributions. In Section~\ref{sec:radial}, we derive the variance-covariance
matrix of any affine equivariant matrix-valued statistic when sampling from a
CES distribution. In Section~\ref{sec:varSCM}, we derive the variance-covariance
matrix of the SCM and provide an application in shrinkage estimation.
Section~\ref{sec:conclusion} concludes. All proofs are kept in the appendix.

\emph{Notation}:
We let  $\I$, $\bo 1$, and $\e_i$ denote the identity matrix, a vector of ones,
and a vector whose $i$th coordinate is one and other coordinates are zero,
respectively.  The notations $(\cdot)^*$, $(\cdot)^\top$, and $(\cdot)^\hop$,
denote the complex conjugate, the transpose, and the conjugate transpose,
respectively. The notations $\Hm{\pdim}$, $\PSDHm{\pdim}$, and $\PDHm{\pdim}$
denote the sets of Hermitian, Hermitian positive semidefinite, and Hermitian
positive definite $p \times p$-dimensional matrices, respectively. We use the
shorthand notation $\var(\bo A) = \var(\ve(\bo A))$ and $\pvar(\bo
A)=\pvar(\ve(\bo A))$ (see Section~\ref{sec:radial} for the definition of
$\pvar$), where $\ve(\bo A) = (\bo a_1^\top \cdots \bo a_\pdim^\top)^\top$ is a
vectorization of $\bo A = (\bo a_1 \cdots \bo a_\pdim)$. When there is a
possibility for confusion, we denote by $\cov_{\bmu,\M}(\cdot,\cdot)$ or
$\E_{\bmu, \M}[\cdot]$ the covariance and expectation of a sample from an
elliptical distribution with mean vector $\bmu$ and covariance matrix $\M$. The
$\pdim^2 \times \pdim^2$ commutation
matrix~\cite{magnusCommutationMatrixProperties1979} is defined by $\commat{p} =
\sum_{i,j} \e_i \e_j^\top \otimes \e_j \e_i^\top$, where $\otimes$ is the
Kronecker product.  The notation $\overset{d}=$ reads ``has the same
distribution as'', $\mathcal U(\fC \mathcal S^{p})$ denotes the uniform
distribution on the complex unit sphere $\fC \mathcal S^{p} = \{\u \in \fC^\pdim
: \|\u\| = 1 \}$ and $\fR_{\geq 0} = \{a \in \fR : a \geq 0\}$.

\section{Complex elliptically symmetric distributions}\label{sec:CES}

A random vector $\x \in \fC^\pdim$ is said to have a circular CES distribution
if and only if it admits the \emph{stochastic representation}
\begin{equation}\label{eq:stochasticrepresentation}
    \x \overset{d}= \bmu + r \bo \M^{1/2} \u,
\end{equation}
where $\bmu = \E[\x]$ is the mean vector, $\M^{1/2} \in \PDHm{\pdim}$ is
the unique Hermitian positive definite square-root of $\M$, $\u \sim \mathcal
U(\fC \mathcal S^{p})$, and $r > 0$ is a positive random variable called the
\emph{modular variate}. Furthermore $r$ and $\u$ are independent. If the
cumulative distribution function of $\x$ is absolutely continuous, the
probability density function exists and is up to a constant of the form
\begin{equation}
    |\M|^{-1} g( (\x - \bmu)^\hop \M^{-1} (\x - \bmu) ),
\end{equation} where $g :
\fR_{\geq 0} \to \fR_{> 0}$ is the \emph{density generator}. We denote this case
by $\x \sim \fC\mathcal E_p(\bmu, \M, g)$. We assume that $\x$ has finite
fourth-order moments, and thus we can assume without any loss of generality that
$\M$ is equal to the covariance matrix
$\var(\x)$~\cite{ollilaComplexEllipticallySymmetric2012a} (implying $\E[r^2] =
p$).  We refer the reader to~\cite{ollilaComplexEllipticallySymmetric2012a} for
a comprehensive account on CES distributions.

The elliptical kurtosis of a CES distribution is defined as
\begin{equation}\label{eq:kappa_complex}
    \kappa = \frac{\E[ r^4]}{p(p+1)} - 1 .
\end{equation}
Elliptical kurtosis shares properties similar to the kurtosis of a circular
complex random variable. Specifically, if $\x \sim \fC \mathcal
N_\pdim(\bmu,\M)$, then $\kappa=0$. This follows by noticing that $ 2 \cdot r^2
\sim \chi^2_{2 \pdim}$, and hence $\E[r^4] = \pdim(\pdim+1)$ and consequently
$\kappa = 0$ in the Gaussian case. The kurtosis of a complex circularly
symmetric random variable $x \in \fC$ is defined as
\begin{equation}\label{eq:kurtosis_complex}
    \kurt(x)= \frac{ \E[|x - \mu |^4 ] }{ (\E[|x - \mu |^2])^2 } -2 ,
\end{equation}
where $\mu = \E[x]$. Similar to the real-valued case, $\kappa$ has a simple
relationship with the (excess) kurtosis
\cite[Lemma~3]{ollilaShrinkingEigenvaluesMestimators2021}:
$
    \kappa = \frac 1 2 \cdot \kurt(x_i),
$
for any $i \in \{1, \ldots,p\}$. We note that the lower bound for the elliptical
kurtosis is $\kappa \geq
-1/(p+1)$~\cite{ollilaComplexEllipticallySymmetric2012a}.

Lastly, we define the \emph{scale} and \emph{sphericity} parameters
\begin{equation}\label{eq:scaleandsphericity}
    \eta = \frac{\tr(\M)}{p} \quad \text{and} \quad \gamma = p
    \frac{\tr(\M^2)}{\tr(\M)^2}.
\end{equation}
The scale is equal to the mean of the eigenvalues. The sphericity measures how
close the covariance matrix is to a scaled identity matrix. The sphericity
parameter gets the value $1$ for the scaled identity matrix and $p$ for a rank
one matrix.

\section{Radial distributions and covariance matrix estimates}\label{sec:radial}

In this section, we derive the variance-covariance matrix of any affine
equivariant matrix-valued statistic. We begin with some definitions. 

The covariance and pseudo-covariance~\cite{neeserProperComplexRandom1993a} of
complex random vectors $\x_1$ and $\x_2$ are defined as
\begin{align*}
    \cov(\x_1,\x_2)
    &= \E \left[ (\x_1 - \E[ \x_1]) (\x_2 - \E[ \x_2])^\hop \right]
    ~\text{and}
    \\
    \pcov(\x_1,\x_2)
    &= \E \left[ (\x_1 - \E[ \x_1]) (\x_2 - \E[ \x_2])^\top \right],
\end{align*}
and together they provide a complete second-order description of associations
between $\x_1$ and $\x_2$. Then $\var(\x)=\cov(\x,\x)$ and $\pvar(\x)=\pcov(
\x,\x)$ are called the covariance matrix and the \emph{pseudo-covariance
matrix} \cite{neeserProperComplexRandom1993a} of $\x$.

A random Hermitian ($\bo A^\hop = \bo A$) matrix $\bo A \in \Hm{p}$ is said to
have a \emph{radial distribution} if $ \bo A \overset{d}= \bo Q \bo A \bo
Q^\hop$ for all unitary matrices $\bo Q$ (so $\bo Q^\hop \bo Q = \bo I$).  The
following result extends the result of~\cite{tylerRadialEstimatesTest1982a} to
the complex-valued case.
\begin{theorem}\label{thm:radial_complex}
    Let a random matrix $\bo A = (a_{ij}) \in \Hm{p}$ have a radial distribution
    with finite second-order moments. Then, there exist real-valued constants
    $\sigma,\tau_1$ and $\tau_2$ with $\tau_1 \geq 0$ and $\tau_2 \geq
    -\tau_1/p$ such that $\E[\bo A]=\sigma \bo I$ with $\sigma=\E[a_{ii}]$ and
    \begin{align}
        \var(\bo A) &= \tau_1 \bo I + \tau_2 \, \ve(\bo I)\ve( \bo I)^\top
        \label{eq:covar} , \\
        \pvar(\bo A) &= \tau_1 \commat{p}+ \tau_2 \, \ve(\bo I)\ve( \bo I)^\top,
        \label{eq:pcovar}
    \end{align}
    where $\tau_1=\var(a_{ij})=\pcov(a_{ij},a_{ji})$ and
    $\tau_2=\cov(a_{ii},a_{jj})=\pcov(a_{ii},a_{jj})$ for all $1\leq i\neq j
    \leq p$.
\end{theorem}

A statistic $\hat \M = \hat \M(\bo X) \in \Hm{p}$ based on an $n \times p$ data
matrix $\bo X = \begin{pmatrix} \x_1 & \cdots & \x_n \end{pmatrix}^\top$ of $n
\geq 1$ observations on $p$ complex-valued variables is said to be \emph{affine
equivariant} if
\begin{equation}\label{eq:affine_eq_complex}
    \hat \M( \bo X \bo A^\top + \bo 1 \bo a^\top)= \bo A \hat \M(\bo X) \bo
    A^\hop
\end{equation}
holds for all $\bo A \in \fC^{p \times p}$ and $\a \in \fC^p$. Suppose that
$\x_1, \ldots, \x_n \subset \fC^\pdim$ is a random sample from a CES
distribution $\fC\mathcal E_p(\bmu,\M,g)$ and that $\hat \M =( \hat \sigma_{ij})
\in \Hm{p}$ is an affine equivariant statistic. Then $\hat \M$ has a stochastic
decomposition
\begin{equation}\label{eq:hatM_radial}
    \hat \M(\bo X)\overset{d}= \M^{1/2} \cdot \hat \M(\bo Z) \cdot \M^{1/2},
\end{equation}
where $\hat \M(\bo Z)$ denotes the value of $\hat \M$ based on a random sample
$\z_1, \ldots, \z_n \subset \fC^\pdim$ from a spherical distribution $\fC
\mathcal E_p(\bo 0,\bo I,g)$. This follows by writing $\bo X \overset{d}= \bo Z
{(\M^{1/2})}^\top + \bom 1 \bmu^\top$ using~\eqref{eq:stochasticrepresentation}
(where $\bo z_i = r_i \bo u_i$) and then applying~\eqref{eq:affine_eq_complex}.
Affine equivariance together with the fact that $\z_i \overset{d} = \bo Q \z_i$
for all unitary matrices $\bo Q$ indicate that $\hat \M(\bo Z)$ has a radial
distribution. This leads to Theorem \ref{thm:cov_affeq} stated below.

\begin{theorem}\label{thm:cov_affeq}
    Let $\hat \M =( \hat \sigma_{ij}) \in \Hm{p}$ be an affine equivariant
    statistic with finite second-order moments, and based on a random sample
    $\x_1, \ldots, \x_n \subset \fC^\pdim$ from a CES distribution $\fC\mathcal
    E_p(\bmu,\M,g)$. Then $\E[\hat \M] = \sigma \M$ with $\sigma = \E_{\bo 0,\bo
    I}[ \hat \sigma_{11}]$, and
    \begin{align}
        \var (\hat \M) &= \tau_1 (\M^* \otimes \M) + \tau_2 \ve(\M)\ve(\M)^\hop,
        \label{eq:var_hatsigma_c} \\
        \pvar(\hat \M) &= \tau_1 (\M^* \otimes \M) \commat{p} + \tau_2
        \ve(\M)\ve(\M)^\top , \label{eq:pvar_hatsigma_c}
    \end{align}
    where $\tau_1=\var_{\bo 0,\bo I}(\hat \sigma_{12})$ and $\tau_2=\cov_{\bo
    0,\bo I}(\hat \sigma_{11}, \hat \sigma_{22}) \geq -\tau_1/p$.
\end{theorem}

There are many statistics for which this theorem applies. Naturally, a prominent
example is the SCM, which we examine in detail in the next section. Other
examples are the complex $M$-estimates of scatter discussed
in~\cite{ollilaComplexEllipticallySymmetric2012a} or the weighted sample
covariance matrices
$
    \bo R = \frac{1}{n} \sum_{i=1}^n u(d_i) (\x_i - \bar \x)(\x_i - \bar
    \x)^\hop,
$
where $d_i = (\x_i - \bar \x)^\hop \bo S^{-1} (\x_i - \bar \x)$ and $u:\fR_{\geq
0} \to \fR_{\geq 0}$. In the special case, when $u(s) = s$, we obtain the fourth
moment matrix used in FOBI~\cite{cardosoSourceSeparationUsing1989} for blind
source separation and in Invariant Coordinate Selection
(ICS)~\cite{tylerInvariantCoordinateSelection2009a}.

\section{Variance-covariance of the SCM and an example in shrinkage
estimation}\label{sec:varSCM}
We now use Theorem \ref{thm:cov_affeq} to derive the covariance matrix and the
pseudo-covariance matrix as well as the MSE of the SCM when sampling from a CES
distribution. This result extends~\cite[Theorem 2 and Lemma
1]{ollilaOptimalShrinkageCovariance2019a} to the complex case.
\begin{theorem}\label{thm:SCM}
    Let the SCM $\mathbf{S} =(s_{ij})$ be computed on an i.i.d. random sample
    $\x_1,\ldots, \x_n \subset \fC^p$ from a CES distribution $\fC \mathcal
    E_p(\bmu,\M,g)$ with finite fourth-order moments and covariance matrix $\M =
    \var(\x_i)$. Then, the covariance matrix and pseudo-covariance matrix of
    $\S$ are as stated in \eqref{eq:var_hatsigma_c} and
    \eqref{eq:pvar_hatsigma_c} with
\begin{align*}
    \tau_1 &= \var_{\bo 0,\bo I}(s_{12}) = \frac{1}{n-1} + \frac{\kappa}{n} , \\
    \tau_2 &= \cov_{\bo 0,\bo I}(s_{11}, s_{22}) = \frac{\kappa}{n},
\end{align*}
where $\ka$ is the elliptical kurtosis in \eqref{eq:kappa_complex}. The MSE is
given by
\begin{align*}
    \MSE(\S) &= \E[\|\S - \M\|_\Fro^2] =
    \Big( \frac{1}{\ndim-1} + \frac{\ka}{\ndim}\Big) \tr(\M)^2 +
    \frac{\ka}{\ndim} \tr(\M^2),
\end{align*}
and the normalized MSE is
    \begin{equation}\label{eq:NMSE}
    \NMSE(\S) = \frac{\MSE(\S)}{\| \M \|_{\Fro}^2} =
    \frac{p}{\gamma}\Big( \frac{1}{\ndim-1} + \frac{\ka}{\ndim}\Big)+
    \frac{\ka}{\ndim},
\end{equation}
where $\gamma$ is the sphericity parameter defined in
    \eqref{eq:scaleandsphericity}.
\end{theorem}

The finite sample performance of the SCM can often be improved by using
shrinkage covariance matrix estimators, which is an active area of research, see
e.g.,~\cite{ledoitWellconditionedEstimatorLargedimensional2004,
ledoitHoneyShrunkSample2004a, duFullyAutomaticComputation2010a,
chenShrinkageAlgorithmsMMSE2010a, ollilaOptimalShrinkageCovariance2019a,
bunRotationalInvariantEstimator2016, ledoitNonlinearShrinkageEstimation2012}.
Consider the simple shrinkage covariance matrix estimation problem,
\begin{equation*}
    \beta_o = \arg \min_{\beta \in \fR} ~ \E[\|\beta \S - \M\|_\Fro^2].
\end{equation*}
Since the problem is convex, we can find $\beta_o$ as solution of
$\frac{\partial}{\partial \beta} \E [ \|\beta \S - \M \|_\Fro^2 ] = 2 \beta
\E[\tr(\S^2)] - 2 \| \M \|_{\Fro}^2 = 0$ which yields
\begin{align}\label{eq:beta_opt}
        \beta_o &=
        \frac{ \| \M \|_{\Fro}^2}{\MSE(\S) + \| \M \|_{\Fro}^2}
        = \frac{1}{\NMSE(\S) + 1},
\end{align}
where we used $\MSE(\S) = \E[\|\S - \M\|_\Fro^2] = \E[\tr(\S^2)] - \| \M
\|_{\Fro}^2$. From~\eqref{eq:beta_opt}, it follows that the optimal scaling term
$\beta_o$ is always smaller than 1 since $\MSE(\S)>0$. Furthermore, $\beta_o$ is
a function of $\gamma$ and $\kappa$ via~\eqref{eq:NMSE}. In the Gaussian case
($\kappa=0$), we obtain $\beta_o = (n-1)/(n-1+ p/\gamma)$.
Figure~\ref{fig:betao} illustrates the effect of $\kappa$ on $\beta_o$ when
$\gamma = 2$, and $n=p=10$. Next we show that the oracle estimator $\S_o =
\beta_o \S$ is uniformly more efficient than $\S$, i.e., $\MSE(\S_0 ) <
\MSE(\S)$ for any $\M \in \PDHm{p}$. First write
\begin{equation}\label{eq:optimal_beta_apu}
     \E \big[ \| \beta_o \S - \M \|^2_{\Fro} \big] = \beta^2_o
     \MSE( \S) + (1-\beta_o)^2 \| \M \|_{\Fro}^2.
\end{equation}
Then from \eqref{eq:beta_opt} we notice that $1- \beta_o = \beta_o \NMSE(\S)$.
Subsituting this into \eqref{eq:optimal_beta_apu} we get
\begin{align*}
    \MSE(\S_o)
    &= \beta_o^2 \MSE(\S) + \beta_o^2 \NMSE(\S)^2 \| \M \|_{\Fro}^{2} \\
    &= \beta_o^2 \MSE(\S) \big(1+ \NMSE(\S) \big)
     = \beta_o \MSE(\S),
\end{align*}
where the last identity follows from fact that $1/\beta_o = 1 + \NMSE(\S) $ due
to~\eqref{eq:beta_opt}. Since $\beta_o <1$ for all $\M \in \PDHm{p}$, it follows
that $\S_o$ is more efficient than $\S$. Efficiency in the case when $\gamma$
and $\kappa$, and hence $\beta_o$ need to be estimated, remains (to the best of
our knowledge) an open problem. However, for certain related shrinkage
estimators the shrinkage intensity can be consistently estimated,
e.g.,~\cite{ledoitWellconditionedEstimatorLargedimensional2004,
ledoitHoneyShrunkSample2004a}.

\newcommand{\nsamples}{10}
\newcommand{\pdimension}{10}
\newcommand{\gamval}{2}
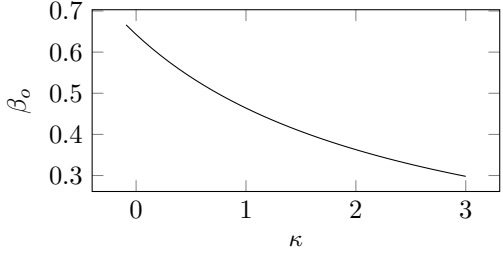
\begin{figure}
    \centering
    \begin{tikzpicture}
        \begin{axis}[
                height = 4.0cm,
                width  = 7cm,
                xlabel = {{\small $\kappa$}},
                ylabel = {$\beta_o$}]
            \addplot[
                domain=(-1/(\pdimension+1)):3, samples=40]
            {1/(1 + \pdimension/\gamval * (1/(\nsamples-1) + x/\nsamples) + x/\nsamples)};
        \end{axis}
    \end{tikzpicture}
    \vspace{-0.1cm}
    \caption{The value of $\beta_o$ as a function of $\kappa$, when $n=p=10$ and
    $\gamma = 2$.}\label{fig:betao}
\end{figure}

In the univariate case ($p=1$),  $\M$ is equal to the variance
$\sigma^2 = \var(x)>0$ of the random variable $x \in \fC$ and the SCM reduces to
the sample variance defined by
$
    s^2 = \frac{1}{\ndim-1} \sum_{i=1}^\ndim |x_i - \bar x|^2.
$
In this case, $\gamma=1$, and
$\beta_o$
in~\eqref{eq:beta_opt} is
\begin{align*}
    \beta_0 &= \frac{n(n-1)}{ \kurt(x) (n-1) +n^2}.
\end{align*}
A similar result was noticed in~\cite{searlsNoteEstimatorVariance1990} for the
real-valued case. If the data is from a complex normal distribution $ \fC
\mathcal N(\mu,\sigma^2$), then $\kurt(x)=0$, and $\beta_o = (n-1)/n$, and hence
$s_o^2 = \beta_o s^2 = \frac{1}{\ndim} \sum_{i=1}^\ndim |x_i - \bar x|^2$, which
equals the Maximum Likelihood Estimate (MLE) of $\sigma^2$. In the real case,
the optimal scaling constant is $\beta_o= (n-1)/(n+1)$ for Gaussian
samples~\cite{goodmanSimpleMethodImproving1953}. Note that when the kurtosis is
large and positive and $n$ is small, then $\beta_o$ can be substantially less
than one and the gain of using $s_o$ can be significant.

\section{Conclusion}\label{sec:conclusion}
We derived the form of the variance-covariance matrix for any affine equivariant
matrix-valued statistics under sampling from CES distributions. We used this
result to derive the variance-covariance matrix and the MSE of the SCM when
finite fourth-order moments exist. An illustrative example in the context of
shrinkage covariance matrix estimation was presented.

\appendix
\subsection{Proof of Theorem~\ref{thm:radial_complex}}\label{app:radial_complex}
The proof follows the same lines as the proof
in~\cite{tylerRadialEstimatesTest1982a} for the real-valued case.
Since $\bo A$ has a radial distribution $\E[\bo A] = \E[a_{ii}] \bo I = \sigma
\bo I$ is obvious. For any unitary matrix $\bo Q = (\bo q_1 \cdots \bo q_p)$, we
have
\begin{align*}
    \var(\bo A)
    &= \var(\ve(\bo A))
    = \var(\ve(\bo Q \bo A \bo Q^\hop))
    \\
    &= \var((\bo Q^* \otimes \bo Q)\ve(\bo A))
    \\
    &=
    (\bo Q^* \otimes \bo Q) \var(\ve(\bo A)) (\bo Q^\top \otimes \bo Q^\hop).
\end{align*}
Let $\{\bo e_i\bo e_j^\top \otimes \bo e_k \bo e_l^\top\}$ be a basis for the
set of $p^2\times p^2$ matrices. Then
\begin{align*}
    \var(\bo A)
    = \sum_{i,j,k,l} \tau_{ijkl} \bo e_i\bo e_j^\top \otimes \bo e_k \bo e_l^\top
    =
    \sum_{i,j,k,l} \tau_{ijkl} \q_i^* \q_j^\top \otimes \q_k \q_l^\hop,
\end{align*}
where $\tau_{ijkl} = \var(a_{ki},a_{lj})$. By choosing $\q_m = \im \bo e_m$
(where $\im$ is the imaginary unit) and $\q_r = \bo e_r$ for some $m \neq r$, we
must have $\tau_{ijkl} = 0$ unless $i=j=k=l$, $i=j \neq k=l$, or $i=k \neq j=l$.
Denote $\tau_0 = \tau_{iiii} = \var(a_{ii})$, $\tau_1 = \tau_{iijj} =
\var(a_{ji})$ and $\tau_2 = \tau_{ijij} = \cov(a_{ii},a_{jj})$. Then
\begin{align*}
    \var(\bo A)
    &=
    \textstyle
    \sum_{i,j} \tau_{1} \bo e_i \bo e_i^\top \otimes \bo e_j \bo e_j^\top
    +
    \sum_{i,j} \tau_{2} \bo e_i \bo e_j^\top \otimes \bo e_i \bo e_j^\top
    \\&\quad
    \textstyle
    + (\tau_0 - \tau_1 - \tau_2) \sum_{i}\bo e_i \bo e_i^\top \otimes \bo e_i
    \bo e_i^\top.
\end{align*}
Note that $\sum_{i,j} \bo e_i \bo e_i^\top \otimes \bo e_j \bo e_j^\top = \bo I$
and $\sum_{i,j} \bo e_i \bo e_j^\top \otimes \bo e_i \bo e_j^\top = \ve(\bo
I)\ve(\bo I)^\top$. Furthermore,
\begin{align*}
    \textstyle
    (\bo Q^* \otimes \bo Q) \sum_{i,j} \bo e_i \bo e_i^\top \otimes \bo e_j \bo
    e_j^\top (\bo Q^\top \otimes \bo Q^\hop)
    &= \bo I
    \\
    \textstyle
    (\bo Q^* \otimes \bo Q) \sum_{i,j} \bo e_i \bo e_j^\top \otimes \bo e_i \bo
    e_j^\top (\bo Q^\top \otimes \bo Q^\hop)
    &= \ve(\bo I)\ve(\bo I)^\top
    \\
    \textstyle
    (\bo Q^* \otimes \bo Q) \sum_{i} \bo e_i \bo e_i^\top \otimes \bo e_i \bo
    e_i^\top (\bo Q^\top \otimes \bo Q^\hop)
    &\neq \bo e_i \bo e_j^\top \otimes \bo e_i \bo e_j^\top.
\end{align*}
From the last inequality, we must have $\tau_0 - \tau_1 - \tau_2 = 0$ and
$\var(\bo A) = \tau_1 \bo I + \tau_2 \ve(\bo I)\ve(\bo I)^\top$ follows.

Regarding the pseudo-covariance, for any unitary $\bo Q$,
\begin{align*}
    \pvar(\bo A) &=
    (\bo Q^* \otimes \bo Q) \pvar(\ve(\bo A)) (\bo Q^\hop \otimes \bo Q^\top),
\end{align*}
which implies
\begin{align*}
    \pvar(\bo A)
    = \sum_{i,j,k,l} \tau'_{ijkl} \bo e_i\bo e_j^\top \otimes \bo e_k \bo e_l^\top
    =
    \sum_{i,j,k,l} \tau'_{ijkl} \q_i^* \q_j^\hop \otimes \q_k \q_l^\top,
\end{align*}
where $\tau'_{ijkl} = \pcov(a_{ki},a_{lj})$. By choosing $\q_m = \im \bo e_m$
and $\q_r = \bo e_r$ for some $m \neq r$, we must have $\tau'_{ijkl} = 0$ except
when $i=j=k=l$, $i=k \neq j=l$, or $i=l \neq j=k$. Let $\tau'_0 = \tau'_{iiii} =
\pvar(a_{ii})$, $\tau'_1 = \tau'_{ijji} = \pcov(a_{ij},a_{ji})$ and $\tau'_2 =
\tau'_{ijij} = \pcov(a_{ii},a_{jj})$. Then,
\begin{align*}
    \pvar(\bo A)
    &=
    \textstyle
    \sum_{i,j} \tau'_{1} \bo e_i \bo e_j^\top \otimes \bo e_j \bo e_i^\top
    +
    \sum_{i,j} \tau'_{2} \bo e_i \bo e_j^\top \otimes \bo e_i \bo e_j^\top
    \\&\quad
    \textstyle
    + (\tau'_0 - \tau'_1 - \tau'_2) \sum_{i}\bo e_i \bo e_i^\top \otimes \bo e_i
    \bo e_i^\top
    \\
    &= \tau'_1 \commat{p} + \tau'_2 \ve(\bo I)\ve(\bo I)^\top
\end{align*}
by similar arguments as with $\var(\bo A)$. Then note that, $\tau_1 =
\var(a_{ji}) = \pcov(a_{ij},a_{ji}) = \tau'_1 \geq 0$ and $\tau_2 =
\cov(a_{ii},a_{jj}) = \pcov(a_{ii},a_{jj}) = \tau'_2$. Lastly, since $\var(\bo
A)$ is positive semidefinite and $\tau_1 \geq 0$,
$
    |\var(\bo A)|
    = |\tau_1 \bo I + \tau_2 \ve(\bo I)\ve(\bo I)^\top|
    = (\tau_1 + \tau_2 p) \tau_1^{p^2-1} \geq 0
$
implies $\tau_2 \geq - \tau_1/p$.
\qed

\subsection{Proof of Theorem~\ref{thm:cov_affeq}}\label{app:cov_affeq}

Since $\hat \M(\bo X)$ is affine equivariant, from~\eqref{eq:hatM_radial}
we have
\begin{align*}
    &\var(\hat \M(\bo X))
    =
    \var( \M^{1/2} \hat \M(\bo Z) \M^{1/2}) \\ 
    &=
    ((\M^{1/2})^* \otimes \M^{1/2})
    \var(\hat \M(\bo Z))
    ((\M^{1/2})^* \otimes \M^{1/2}).
\end{align*}
From Theorem~\ref{thm:radial_complex}, $\var(\hat \M(\bo Z))$ is of the form
\eqref{eq:covar}. Since
\begin{align*}
    &((\M^{1/2})^* \otimes \M^{1/2})
    \bo I
    ((\M^{1/2})^* \otimes \M^{1/2}) = (\M^* \otimes \M) ~\text{and}
    \\
    &((\M^{1/2})^* \otimes \M^{1/2})
    \ve(\bo I)\ve(\bo I)^\top
    ((\M^{1/2})^* \otimes \M^{1/2})
    \\&= \ve(\M)\ve(\M)^\hop,
\end{align*}
we obtain~\eqref{eq:var_hatsigma_c}. Similarly,
\begin{align*}
    &\pvar( \hat \M(\bo X)) = \pvar(\M^{1/2} \hat \M(\bo Z) \M^{1/2})
    \\
    &=
    ((\M^{1/2})^* \otimes \M^{1/2})
    \pvar(\hat \M(\bo Z))
    (\M^{1/2} \otimes (\M^{1/2})^*),
\end{align*}
where $\pvar(\hat \M(\bo Z))$ is of the form~\eqref{eq:pcovar}. Since
\begin{align*}
    &((\M^{1/2})^* \otimes \M^{1/2})
    \commat{p}
    (\M^{1/2} \otimes (\M^{1/2})^*)
    = (\M^* \otimes \M) \commat{p},
    \\
    &((\M^{1/2})^* \otimes \M^{1/2})
    \ve(\bo I)\ve(\bo I)^\top
    (\M^{1/2} \otimes (\M^{1/2})^*)
    \\&= \ve(\M)\ve(\M)^\top,
\end{align*}
we obtain~\eqref{eq:pvar_hatsigma_c}.
\qed

\subsection{Proof of Theorem 3}\label{app:complexSCM}

The proof is similar to the proof
of~\cite[Theorem~2]{ollilaOptimalShrinkageCovariance2019a} that was derived for
the real-valued case. Write the SCM as $\S = (s_{ij})=(n-1)^{-1} \bo
X^\top \bo H \bo X^{*}$, where $\bo H = \bo I - \frac{1}{n}\bo 1 \bo 1^\top $ is
the centering matrix.  Write $\a=\bo X \bo e_q$ and $\bo b=\bo X \bo e_r$ for $q
\neq r$. Then note that $s_{qr} = (n-1)^{-1} \a^\top \bo H \bo b^{*}$.  Hence,
\begin{align}\label{eq:tau1_proof_c1}
    \tau_1 &= \var(s_{qr}) = \var((n-1)^{-1} \a^\top \bo H \bo b^*) \notag
    \\
    &=
    (n-1)^{-2} \var(\a^\top \bo H \bo b^*).
\end{align}
Then note that
\begin{align}
    \var(\a^\top \bo H \bo b^*)
    &= \var(\tr( \bo H \bo b^* \a^\top ) ) \\
    &= \var( \veco{\bo H}^\top\veco{\bo b^* \a^\top}) \notag \\
    &= \veco{\bo H}^\top \var \big( \veco{\bo b^* \a^\top} \big) \veco{\bo H}.
    \label{eq:tau1_proof_c2}
\end{align}
Since for $\x_i \sim\fC \mathcal{E}_p(\bo 0,\bo I,g)$ we have $\x_i
\overset{d}= r_i \u_i$, where $r_i \overset{d}= \|\x_i\|$ is independent of
$\u_i = {(u_{i1},\ldots,u_{ip})}^\top \sim \mathcal U(\fC \mathcal S^{p})$, we
can write $\a=\bo X \bo e_q = {(r_1 u_{1q}, r_2 u_{2q}, \ldots, r_n u_{nq})}^\top$,
and similarly for $\bo b$. The $kl$th element of the $ij$th block (i.e., the
$ijkl$th element) of the $n^2 \times n^2$ matrix $\var\big(\veco{\bo b^*
\a^\top} \big)$ is
\begin{align*}
    \cov\left( b_k^* a_i , b_l^* a_j\right)
    &= \E\left[ r_{k} r_{i} r_{l} r_{j} u_{kr}^*u_{iq} u_{lr} u_{jq}^*
    \right],
\end{align*}
where we used that $\E[ \bo b^* \a^\top]=\bo 0$. Then note that
\begin{align*}
    &\E\left[ |u_{iq}|^2| u_{ir}|^2\right] = \frac{1}{p(p+1)}, \\
    &\E\left[| u_{iq}|^2\right]            = \frac{1}{p}, ~\text{and}~
        \E\left[ |u_{iq}|^4\right] = \frac{2}{p(p+1)},
\end{align*}
while all other moments up to fourth-order vanish. This and the fact that
$\E\left[r_i^4\right] = (1+\kappa) p(p+1)$ due to~\eqref{eq:kappa_complex},
implies that the only non-zero elements of $\var\big(\veco{\bo b^* \a^\top}
\big)$ are
\begin{align*}
    \E[r_i^4] \E[ |u_{ir}|^2 |u_{iq}|^2] &= 1+\kappa
        &\text{for}~i&=j=k=l,~\text{and}
    \\
    \E[r_{i}^2]\E[r_{k}^2] \E[u_{ir}^2]\E[u_{kq}^2] &= 1
        &\text{for}~i&=j\neq k=l,
\end{align*}
and hence
\begin{equation}\label{eq:tau1_proof_c3}
    \var\big(\veco{\bo b^* \a^\top} \big) 	= \bo I + \kappa \sum_{i=1}^n \bo
    e_i \bo e_i^\top \otimes \bo e_i \bo e_i^\top.
\end{equation}
This together with \eqref{eq:tau1_proof_c1} and \eqref{eq:tau1_proof_c2} yields
\begin{align*}
    \tau_1 &= \frac{1}{(n-1)^{2}} \veco{\bo H}^\top \left(	\bo I + \kappa
    \sum_{i=1}^n \bo e_i \bo e_i^\top \otimes \bo e_i \bo
    e_i^\top\right)\veco{\bo H} \\
    &= \frac{1}{n-1} + \frac{\kappa}{n},
\end{align*}
where we used $\veco{\bo H}^\top \veco{\bo H} = n - 1$ and
\begin{align*}
    \sum_{i=1}^n \veco{\bo H}^\top (\bo e_i \bo e_i^\top \otimes \bo e_i \bo e_i^\top)
    \veco{\bo H} = \sum_{i=1}^n h_{ii}^2 = \frac{{(n-1)}^2}{n}.
\end{align*}

Next, we find the expression for $\tau_2=\cov(s_{qq},s_{rr}) = (n-1)^{-2}
\cov(\a^\top \bo H \a^*, \bo b^\top \bo H \bo b^*)$. Since $\E[s_{qq}] =
\E[(n-1)^{-1} \a^\top \bo H \a^*] = 1$ for any $q$,
\begin{align*}
    \tau_2 &=
    (n-1)^{-2}\E[ \a^\top \bo H \a^* \bo b^\top \bo H \bo b^* ] - 1 \\
    &= (n-1)^{-2} \tr \big( \E \big[ \, \bo H (\bo b^* \a^\top)^\hop \bo H
    (\bo b^* \a^\top) \, \big] \big) - 1 \\
    &= (n-1)^{-2} \tr \Big( (\bo H \otimes \bo H) \E \Big[ \veco{\bo b^*
    \a^\top} \veco{\bo b^* \a^\top}^\hop \Big] \Big) - 1.
\end{align*}
The expression in the expectation is equal to~\eqref{eq:tau1_proof_c3}, and so
\[
    \tau_2 = \frac{1}{(n-1)^2} \tr \Big( (\bo H \otimes \bo H) \Big( \bo I +
    \kappa \sum_{i=1}^n \bo e_i \bo e_i^\top \otimes \bo e_i \bo e_i^\top \Big)
    \Big) - 1 =
    \frac{\kappa}{n},
\]
where we used that $\tr(\bo H \otimes \bo H) = \tr(\bo H)^2=(n-1)^2$ and
\begin{align*}
    \sum_{i=1}^n\tr\left((\bo H \otimes \bo H)(\bo e_i \bo e_i^\top \otimes \bo e_i\bo e_i^\top
    )\right)=\sum_{i=1}^n h_{ii}^2.
\end{align*}
This completes the proof for $\tau_1$ and $\tau_2$. By
Theorem~\ref{thm:cov_affeq}, we have $\var (\S) = \tau_1 (\M^* \otimes \M) +
\tau_2 \ve(\M)\ve(\M)^\hop$, and hence
\begin{align*}
    \MSE(\S) &=\tr(\var(\S)) \\ 
&= \tr \left(\tau_1 (\M^* \otimes \M) + \tau_2 \ve(\M)\ve(\M)^\hop\right)
    \\
    &= \tau_1 \tr(\M)^2 + \tau_2 \tr(\M^2), 
\end{align*}
where the last identity follows from using $\tr(\M^* \otimes \M) = \tr(\M)^2$.
This gives the stated expression for the MSE. \qed

\clearpage
\bibliographystyle{IEEEtran}
\input{complexSCM.bbl}

\end{document}